\documentclass[12pt]{article}
\usepackage{amssymb}
\usepackage{amsfonts}
\usepackage{amsmath}
\usepackage[dvips,bottom=1.4in,right=1in,top=1in, left=1in,nohead]{geometry}
\usepackage{setspace}
\usepackage{soul}

\newtheorem{theorem}{Theorem}[section]

\newtheorem{definition}[theorem]{Definition}

\newtheorem{proposition}[theorem]{Proposition}
\newtheorem{remark}[theorem]{Remark}

\newenvironment{proof}[1][Proof]{\noindent\textbf{#1.} }{\ \rule{0.5em}{0.5em}}
\numberwithin{equation}{section}
\input{tcilatex}

\begin{document}

\title{\textbf{Sharp estimates for the global attractor of scalar
reaction-diffusion equations with a Wentzell boundary condition}}
\author{Ciprian G. Gal \\
Department of Mathematics\\
University of Missouri,\\
Columbia, MO 65211, USA\\
ciprian@math.missouri.edu}
\maketitle

\begin{abstract}
\noindent In this paper, we derive optimal upper and lower bounds on the
dimension of the attractor $\mathcal{A}_{W}$\ for scalar reaction-diffusion
equations with a Wentzell (dynamic) boundary condition. We are also
interested in obtaining explicit bounds on the constants involved in our
asymptotic estimates, and to compare these bounds to previously known
estimates for the dimension of the global attractor $\mathcal{A}_{K},$ $K\in
\left\{ D,N,P\right\} $, of reaction-diffusion equations subject to
Dirichlet, Neumann and periodic boundary conditions. The explicit estimates
we obtain show that the dimension of the global attractor \ $\mathcal{A}_{W}$
is of different order than the dimension of $\mathcal{A}_{K},$ for each $%
K\in \left\{ D,N,P\right\} ,$ in all space dimensions that are greater or
equal than three.
\end{abstract}

\doublespacing

\footnotetext[1]{%
2000 Mathematics Subject Classification: 35K57, 35K55, 35B40, 35B41, 37L30,
35Q80, 35Q86}

\section{Introduction}


It is well-known that the long-time behaviour of solutions of partial
differential equations arising in mathematical physics can, in many cases,
be described in terms of global attractors of the associated semigroups (see 
\cite{BV, CV, Ha, T} and references therein). For a large class of equations
of mathematical physics, including parabolic partial differential equations
modelling reaction, diffusion and drift, hyperbolic type equations, and so
on, the corresponding attractor has finite Hausdorff and fractal dimensions.
Thus, the dynamics on the attractor happens to be finite-dimensional, even
though the system is governed by a set of partial differential equations. As
the dimension of the attractor is indicative of the number of degrees of
freedom needed to simulate a given dynamical system, it is then crucial to
obtain more realistic estimates for its dimension in terms of observable
physical quantities.

Aside from some applied motivation, much of the mathematical interest
nowadays is centered on the dynamics of boundary value problems with \emph{%
static }boundary conditions of Dirichlet and Neumann-Robin type, or even
periodic boundary conditions. The influence of these dissipative boundary
conditions on a given model has only been recently investigated in
connection with a class of reaction-diffusion systems. In \cite{HR}, a first
contribution is made to the understanding of this problem with a Robin
boundary condition. In particular, it is shown, for a fixed nonlinearity,
how the flow defined by the reaction-diffusion system depends on the
interaction between diffusion $\nu $\ and another parameter $\theta $\
involved in the boundary condition (cf. also \cite{HR2}). A classification
of points in $\left( \nu ,\theta \right) $-space, as structurally stable, or
bifurcation points, for a one-dimensional scalar reaction-diffusion equation
with a cubic nonlinearity is discussed in detail in \cite{HR}. Other studies
on the influence of boundary conditions upon the solution structures of
partial differential equations have also been done by other scientists.
These studies have analyzed the detailed effect of boundary conditions on
the structure of global attractors (see, e.g., \cite{DMO, HS, MT, SG}). If
the equilibrium is nonhyperbolic and a bifurcation occurs, the structure of
attractors may vary with respect to boundary conditions. This has been
observed in the analysis of pattern formation in a 1D reaction-diffusion
system \cite{DMO}, in lattice systems \cite{S}, in the study of steady state
bifurcations \cite{HS, MT}, and finally in \cite{SG}, on mode-jumping of the
von Karman equations.

Although the global attractors of these systems will depend, for a given
nonlinearity, on the choice of the boundary conditions, their finite
dimension does generally \emph{not}. This result can be easily formulated
for a scalar reaction-diffusion equation, as follows. Consider the parabolic
partial differential equation%
\begin{equation}
\partial _{t}u=\nu \Delta u-f\left( u\right) +\lambda u+g,\text{ }\left(
x,t\right) \in \Omega \times \left( 0,+\infty \right) ,  \label{she}
\end{equation}%
where $u=u\left( x,t\right) \in \mathbb{R}$, $\Omega \subset \mathbb{R}^{n}$%
, $n\geq 1$, is a bounded domain with sufficiently smooth boundary $\Gamma ,$
$g=g\left( x\right) $, and $\nu $, $\lambda $ are positive constants. The
function $f:\mathbb{R\rightarrow R}$ is assumed to be $C^{1,1}$, that is,
continuous and with a Lipschitz continuous first derivative, which
satisfies, among other natural growth conditions (see, e.g., \cite[Chapter II%
]{CV}),%
\begin{equation*}
f^{^{\prime }}\left( y\right) \geq -c_{f},\text{ for all }y\in \mathbb{R}%
\text{, for some }c_{f}>0.
\end{equation*}%
We may ask that $u$ satisfy either a Dirichlet ($K=D$) boundary condition or
a Neumann ($K=N$)\ boundary condition, and even a periodic ($K=P$) boundary
condition. It is well-known that equation (\ref{she}), supplemented with an
appropriate initial condition, generates a semigroup $\left\{ S_{t}\right\} $
acting on a suitable phase space $H$. This semigroup possesses the global
attractor $\mathcal{A}_{K}$, which may depend on the choice of the boundary
conditions, and $\mathcal{A}_{K}$ has finite fractal dimension for each $%
K\in \left\{ D,N,P\right\} $. In particular, the Haussdorf and fractal
dimensions of $\mathcal{A}_{K},$ for any $K\in \left\{ D,N,P\right\} $,
satisfy the following upper and lower bounds:%
\begin{equation}
c_{0}\left( \frac{\lambda }{\nu }\right) ^{n/2}\left\vert \Omega \right\vert
\leq \dim _{H}\mathcal{A}_{K}\leq \dim _{F}\mathcal{A}_{K}\leq c_{1}\left( 1+%
\frac{c_{f}+\lambda }{\nu }\left\vert \Omega \right\vert ^{2/n}\right)
^{n/2},  \label{lb1}
\end{equation}%
for some positive constants $c_{0},c_{1}$ that depend only on $n,$ $f$ and
the shape of $\Omega $ (see, e.g., \cite[Chapter III]{BV}; cf. also \cite{CV}%
, \cite[Chapter VI]{T}). Here, $\left\vert \Omega \right\vert $ stands for
the Lebesgue measure of $\Omega $. For a fixed domain $\Omega $, we observe
that these estimates are sharp with respect to $\nu \rightarrow 0^{+}$ (for
each fixed $\lambda >0$), or sufficiently large $\lambda $ (for each fixed $%
\nu >0$). Hence, these bounds for the dimension of $\mathcal{A}_{K}$\ are of
the same order for each $K\in \left\{ D,N,P\right\} $. These remarkable
estimates\ also depend linearly on the "volume" of the spatial domain $%
\Omega $, which is consistent with physical intuition. This property of the
dimension of the attractor has not been proved for all equations, such as,
the Kuramoto-Sivashinski equation.

Our main goal in this paper is to investigate the dependance of the
dimension of the global attractor for equation (\ref{she}) subject to a
completely new class of boundary conditions, which are sometimes dubbed as 
\emph{Wentzell} boundary conditions, and which have some applications in
probability theory, specifically, Markov processes. But what are they
really? To put them into a context, let $L$ be an elliptic differential
operator of the second order (e.g., $L=\nu \Delta $) with coefficients that
are well-defined over $\overline{\Omega }$. It is known that there is a
one-to-one correspondence between $\left( C_{0}\right) $-semigroups and
Markov processes in $\overline{\Omega }$ which are homogeneous in time and
satisfy the condition of Feller \cite{Fe1} (that is, the range of the
resolvent operator coincides with a prescribed set). Thus, to each such
Markov process there is a corresponding semigroup of operators%
\begin{equation*}
T_{t}v\left( x\right) =\int\limits_{\overline{\Omega }}v\left( y\right)
P\left( t,x,dy\right) ,
\end{equation*}%
where the Markov transition function $P\left( t,x,B\right) $ satisfies $%
P\left( t,x,B\right) \geq 0,$ for $t\geq 0,$ $x\in \overline{\Omega }$ and
any Borel set $B\subseteq \overline{\Omega }$. As a function of $B$, $%
P\left( t,x,\cdot \right) $ is a probability measure. What are the most
general boundary conditions which restrict the given operator $L$\ (more
correctly, its closure) to the infinitesimal operator of a semigroup of
positive contraction operators acting on $C\left( \overline{\Omega }\right) $%
? Wentzell \cite{W} gave a partial answer to this question in higher space
dimensions by finding a sufficiently large class of boundary conditions
which involve differential operators on the boundary that are of the same
order as the operator acting in $\Omega $. He discovered the following form
of boundary conditions:%
\begin{equation}
Lu+\nu b\partial _{\mathbf{n}}^{L}u=0\text{, on }\Gamma \times \left(
0,+\infty \right) ,  \label{wbc}
\end{equation}%
where $\mathbf{n}$ denotes the outward normal at $\Gamma $, $b$ is a
positive constant and $\partial _{\mathbf{n}}^{L}u$ is the outward co-normal
derivative of $u$ with respect to $L$. We refer also to the pioneering work
of \cite{Fe2}, for generation theorems for $L$ with Wentzell boundary
conditions in one space dimension. Until the work of \cite{FGGR}, the study
of the operator $L$\ with Wentzell boundary conditions was usually confined
to generation properties of this operator in the space $C\left( \overline{%
\Omega }\right) $. In 2002, the authors in \cite{FGGR} have found a way to
introduce the Wentzell boundary condition (\ref{wbc}) in an $L^{p}$-context,
which led to the discovery of the natural space for these type of problems
(see Section 2). The reader is referred to \cite{CFGGOR, Gi2} for an
extensive survey of these results and some history.

For the homogeneous linear heat equation (\ref{she}) (that is, $f=g=\lambda
=0$), the Wentzell boundary condition (\ref{wbc})\ is equivalent to a purely
differential equation of the form%
\begin{equation}
\partial _{t}u+\nu b\partial _{\mathbf{n}}u=0,\quad \text{on}\;\Gamma \times
\left( 0,\infty \right) .  \label{dyn}
\end{equation}%
Thus, the main attraction here is that there is a dynamic element introduced
into the boundary condition. The heat equation, supplemented by either (\ref%
{wbc}) or (\ref{dyn}), corresponds to the situation where there is a heat
source (if $b>0$) or sink (if $b<0$) acting on the boundary $\Gamma $.
Mathematically speaking, this kind of conditions (\ref{dyn}) arises due to
the presence of additional boundary terms in the free energy, which must
also account for the action of a source on $\Gamma $\ (see \cite{GW}). We
refer the reader to \cite{Gi} (cf. also \cite{GG0}), for an extensive
derivation and physical interpretation of (\ref{dyn}) for (\ref{she}). For
the nonlinear parabolic equation (\ref{she}), the boundary condition (\ref%
{wbc}) can be formally be transformed into a condition of the form%
\begin{equation}
\partial _{t}u+\nu b\partial _{\mathbf{n}}u+f\left( u\right) -\lambda u=g,%
\text{ on }\Gamma \times \left( 0,\infty \right) .  \label{dyn2}
\end{equation}%
Generally, one may replace $f-\lambda $ in (\ref{dyn2}) by another arbitrary
function $h$, satisfying suitable assumptions. With more sophisticated
arguments, using techniques from semigroup theory, and a variation of
parameter formula, it is possible to prove that the regularity of the
solution for (\ref{she}),(\ref{dyn2}) increases as $f$, $\Omega $ and $g$
become more regular (see Section 2). In particular, for $g=0$, if $\Omega $
is a bounded $\mathcal{C}^{\infty }$ domain and $f$ is a $\mathcal{C}%
^{\infty }$ function, regularity theory implies that the solution $u\left(
t\right) $ to (\ref{she}),(\ref{dyn2}) belongs $H^{k}\left( \Omega \right) ,$
for all $k\geq 0$ and all positive times. At least in this case, the
boundary condition (\ref{wbc}), for equation (\ref{she}), is equivalent to
the boundary condition (\ref{dyn2}). However, in general, this may not be so
if the solution, for the semilinear problem (\ref{she}) and the Wentzell
condition (\ref{dyn2}), is not smooth enough. Since we wish to treat the
most general case, by imposing the least regularity assumptions on $f,$ $g$
and $\Omega ,$ we will devote our attention only to the study of (\ref{she}%
), subject to linear boundary conditions of the form (\ref{dyn}). Our
results below can be immediately extended to other classes of nonlinear
Wentzell boundary conditions (see, e.g., \cite{GW} and references therein).
Boundary conditions of the form (\ref{dyn2}) arise for many known equations
of mathematical physics. They are motivated by heat control problems
formulated in the book of Duvaut and Lions \cite{DL}, problems in
phase-transition phenomena \cite{CGGM, G, GG1, GM, GM2, GMS, GMS2, MZ1, MZ2}
(and their references), special flows in hydrodynamics \cite{FL, GW, QWS, MW}%
, Stefan problems \cite{A, Ma, RSY}, models in climatology \cite{MW2}, and
many others. The reader is referred to \cite{GG_b} for a more complete list
of references involving the application of such boundary conditions to
real-world phenomena.

By keeping our treatment of the boundary condition simple, we wish to prove
that the problem (\ref{she}), (\ref{dyn}) generates a dynamical system on a
suitable phase-space, possessing a finite dimensional global attractor $%
\mathcal{A}_{W}.$ Then, we establish that the Haussdorf and fractal
dimensions of $\mathcal{A}_{W}$ satisfy the following \ upper and lower
bounds:%
\begin{equation}
c_{1}\left( \frac{\lambda }{C_{W}\left( \Omega ,\Gamma \right) \nu }\right)
^{n-1}\leq \dim _{H}\mathcal{A}_{W}\leq \dim _{F}\mathcal{A}_{W}\leq
c_{2}\left( 1+\frac{c_{f}+\lambda }{C_{W}\left( \Omega ,\Gamma \right) \nu }%
\right) ^{n-1},  \label{lb2}
\end{equation}%
for $n\geq 2,$ and%
\begin{equation}
c_{3}\left( \frac{\lambda }{C_{D}\left( \Omega \right) \nu }\right)
^{1/2}\leq \dim _{H}\mathcal{A}_{W}\leq \dim _{F}\mathcal{A}_{W}\leq
c_{4}\left( 1+\frac{c_{f}+\lambda }{C_{D}\left( \Omega \right) \nu }\right)
^{1/2},  \label{lb3}
\end{equation}%
in one space dimension. The positive constants $c_{i},$ $i=1,...,4,$ depend
only on $n,$ $f$ and the shape of $\Omega ,$ while explicit estimates and
formulas for $C_{W}\left( \Omega ,\Gamma \right) $ and $C_{D}\left( \Omega
\right) ,$ respectively, are provided in the Appendix. We note again that,
for a fixed domain $\Omega $, these estimates are sharp with respect to $\nu
\rightarrow 0^{+}$ (for each fixed $\lambda >0$), and for sufficiently large 
$\lambda $ (if $\nu >0$ is fixed). We remark that the bounds we obtain in (%
\ref{lb2})-(\ref{lb3}) are quite simple and their explicit dependance on the
physical parameters is transparent. Moreover, a careful analysis of the
constants involved in (\ref{lb2}) yields the following more explicit
two-sided estimate,%
\begin{equation}
c_{1}^{^{\prime }}\left( \frac{\lambda }{\nu b}\right) ^{n-1}\left\vert
\Gamma \right\vert \leq \dim _{H}\mathcal{A}_{W}\leq \dim _{F}\mathcal{A}%
_{W}\leq c_{2}^{^{\prime }}\left( 1+\frac{c_{f}+\lambda }{\nu b}\left\vert
\Gamma \right\vert ^{1/\left( n-1\right) }\right) ^{n-1},  \label{lb4}
\end{equation}%
in all space dimensions $n\geq 3$. It is worth pointing out that the bounds
in (\ref{lb4}) are proportional to the "surface area" $\left\vert \Gamma
\right\vert $ of $\Gamma ,$ and \emph{not} the "volume" $\left\vert \Omega
\right\vert $\ of $\Omega .$ This is remarkable; most nonlinear equations
arising in mathematical physics, involving the Laplacian on bounded domains,
have the dimension of the attractor of the order of $\left\vert \Omega
\right\vert ^{\alpha },$ for some $\alpha >0$ and for sufficiently large
domains. This property may have profound implications in the prediction of
weather and climate. The reader is referred to Section 4 where this
interesting physical observation is further discussed for the balance
equations governing the large-scale oceanic motion.

Our paper is organized as follows. In Section \ref{ub}, we obtain upper
bounds (cf. Theorem \ref{main1})\ for the fractal dimension of the global
attractor for equation (\ref{she}) with dynamic boundary conditions of the
form (\ref{dyn}). In Section \ref{lbb}, we employ the same technique of \cite%
{BV} to derive a lower bound for the dimension of the unstable manifold of a
constant stationary solution $u^{\ast }$ of (\ref{she}), (\ref{dyn}). As a
consequence, we find a lower bound on the dimension of $\mathcal{A}_{W}$
(see Theorem \ref{lbbb}). Finally, in the Appendix, we recall some useful
results on the so-called Wentzell Laplacian, and prove an auxiliary
inequality, namely, we derive some kind of Sobolev-Lieb-Thirring inequality
that is required to prove the upper bound in (\ref{lb2}).

\section{Upper bounds on the dimension}

\label{ub}

We use the standard notation and facts from the dynamic theory of parabolic
equations (see, for instance, \cite{CFGGOR}, \cite{FGGR}, \cite{GG1}, \cite%
{GW}). We denote by $\left\Vert \cdot \right\Vert _{p}$ and $\left\Vert
\cdot \right\Vert _{p,\Gamma },$ the norms on $L^{p}\left( \Omega \right) $
and $L^{p}\left( \Gamma \right) ,$ respectively. In the case $p=2$, $\langle
\cdot ,\cdot \rangle _{2}$ stands for the usual scalar product. The norms on 
$H^{r}\left( \Omega \right) $ and $H^{r}\left( \Gamma \right) $ are
indicated by $\left\Vert \cdot \right\Vert _{H^{r}\left( \Omega \right) }$
and $\left\Vert \cdot \right\Vert _{H^{r}\left( \Gamma \right) }$,
respectively, for any $r>0$.

The natural phase-space for problem (\ref{she}), (\ref{dyn}) is%
\begin{equation*}
\mathbb{X}^{p}:=L^{p}(\Omega )\oplus L^{p}(\Gamma )=\{F=\binom{f}{g}:f\in
L^{p}(\Omega ),\;g\in L^{p}(\Gamma )\},
\end{equation*}%
for all $p\in \left[ 1,\infty \right] $, endowed with the norm%
\begin{equation}
\left\Vert F\right\Vert _{\mathbb{X}^{p}}^{p}=\int_{\Omega }\left\vert
f\left( x\right) \right\vert ^{p}dx+\int_{\Gamma }\left\vert g(x)\right\vert
^{p}\frac{dS}{b},\text{ }b>0,  \label{xp}
\end{equation}%
if $p\in \lbrack 1,\infty ),$ and 
\begin{equation*}
\Vert F\Vert _{\mathbb{X}^{\infty }}:=\Vert f\Vert _{L^{\infty }(\Omega
)}+\Vert g\Vert _{L^{\infty }(\Gamma )}.
\end{equation*}%
In the definition of $\mathbb{X}^{p}$, $dx$ denotes the Lebesgue measure on $%
\Omega ,$ and $dS$ denotes the natural surface measure on $\Gamma $.
Moreover, we have \cite{FGGR},%
\begin{equation*}
\mathbb{X}^{p}=L^{p}\left( \overline{\Omega },d\mu \right) ,\text{ }p\in %
\left[ 1,\infty \right] ,
\end{equation*}%
where the measure $d\mu =dx_{\mid \Omega }\oplus \frac{dS}{b}_{\mid \Gamma
}, $ on $\overline{\Omega },$ is defined for any measurable set $B\subset 
\overline{\Omega }$ by $\mu (B)=|B\cap \Omega |+\left\vert B\cap \Gamma
\right\vert $. The Dirichlet trace map $\mathit{Tr}_{D}:C^{\infty }\left( 
\overline{\Omega }\right) \rightarrow C^{\infty }\left( \Gamma \right) ,$
defined by $\mathit{Tr}_{D}\left( u\right) =u_{\mid \Gamma }$ extends to a
linear continuous operator $\mathit{Tr}_{D}:H^{r}\left( \Omega \right)
\rightarrow H^{r-1/2}\left( \Gamma \right) ,$ for all $r>1/2$, which is onto
for $1/2<r<3/2.$ This map also possesses a bounded right inverse $\mathit{Tr}%
_{D}^{-1}:H^{r-1/2}\left( \Gamma \right) \rightarrow H^{r}\left( \Omega
\right) $ such that $\mathit{Tr}_{D}\left( \mathit{Tr}_{D}^{-1}\psi \right)
=\psi ,$ for any $\psi \in H^{r-1/2}\left( \Gamma \right) $. Identifying
each function $v\in C\left( \overline{\Omega }\right) $ with the vector $V=%
\binom{v}{\mathit{Tr}_{D}\left( v\right) }\in C\left( \overline{\Omega }%
\right) \times C\left( \Gamma \right) $, it follows that $C\left( \overline{%
\Omega }\right) $ is a dense subspace of $\mathbb{X}^{p},$ for every $p\in
\lbrack 1,\infty ),$ and a closed subspace of $\mathbb{X}^{\infty }.$
Finally, we can also introduce the subspaces of $H^{r}\left( \Omega \right)
\times H^{r-1/2}\left( \Gamma \right) ,$%
\begin{equation*}
\mathbb{V}_{r}:=\left\{ \binom{u}{\psi }\in H^{r}\left( \Omega \right)
\times H^{r-1/2}\left( \Gamma \right) :\mathit{Tr}_{D}\left( u\right) =\psi
\right\} ,
\end{equation*}%
for every $r>1/2,$ and note that we have the following dense and compact
embeddings $\mathbb{V}_{r_{1}}\subset \mathbb{V}_{r_{2}},$ for any $%
r_{1}>r_{2}>1/2$. The linear subspace $\mathbb{V}_{r}$ is densely and
compactly embedded into $\mathbb{X}^{2},$ for any $r>1/2$. We emphasize that 
$\mathbb{V}_{r}$ is not a product space and that, due to the boundedness of
the trace operator $\mathit{Tr}_{D},$ $\mathbb{V}_{r}$ is topologically
isomorphic to $H^{r}\left( \Omega \right) $ in the obvious way.

We begin by stating all the hypotheses on $f$ and $g$ that we need. We
assume that $g\in L^{2}\left( \Omega \right) $ and the following conditions
for $f\in C^{1}\left( \mathbb{R}\text{,}\mathbb{R}\right) $ hold:%
\begin{equation}
f^{^{\prime }}\left( y\right) >-c_{f},\text{ for all }y\in \mathbb{R}\text{,}
\label{n1}
\end{equation}%
\begin{equation}
\eta _{1}\left\vert y\right\vert ^{p}-C_{f}\leq f\left( y\right) y\leq \eta
_{2}\left\vert y\right\vert ^{p}+C_{f},  \label{n2}
\end{equation}%
for some $\eta _{1}$, $\eta _{2}>0,$ $C_{f}\geq 0$ and $p>2$.

We have the following rigorous notion of weak solution to (\ref{she}), (\ref%
{dyn}), with initial condition $u\left( 0\right) =u_{0},$ as in \cite{GW}.

\begin{definition}
\label{weak}The pair $U\left( t\right) =\binom{u\left( t\right) }{\psi
\left( t\right) }$ is said to be a weak solution if $\psi \left( t\right) =%
\mathit{Tr}_{D}\left( u\right) $ for almost all $t\in \left( 0,T\right) ,$
for any $T>0$, and $U$ fulfills%
\begin{align}
U& \in C\left( \left[ 0,T\right] ;\mathbb{X}^{2}\right) \cap L^{\infty
}\left( 0,T;\mathbb{V}_{1}\right) \cap L^{p}\left( \Omega \times \left(
0,T\right) \right) ,  \label{reg} \\
u& \in H_{loc}^{1}(0,\infty ;L^{2}\left( \Omega \right) ),\text{ }\psi \in
H_{loc}^{1}(0,\infty ;L^{2}\left( \Gamma \right) ),  \notag \\
\partial _{t}U& \in L^{2}\left( 0,T;\mathbb{V}_{1}^{\ast }\right) ,  \notag
\end{align}%
such that the identity%
\begin{equation}
\left\langle \partial _{t}U,\Xi \right\rangle _{\mathbb{X}^{2}}+\nu
\left\langle \nabla u,\nabla \sigma \right\rangle _{2}+\left\langle f\left(
u\right) -\lambda u,\sigma \right\rangle _{2}=\left\langle g,\sigma
\right\rangle _{2},  \notag
\end{equation}%
holds almost everywhere in $\left( 0,T\right) $, for all $\Xi =\binom{\sigma 
}{\varpi }\in \mathbb{V}_{1}$. Moreover, we have, in the space $\mathbb{X}%
^{2}$,%
\begin{equation}
U\left( 0\right) =\binom{u_{0}}{v_{0}}=:U_{0},  \label{ini}
\end{equation}%
where $u\left( 0\right) =u_{0}$ almost everywhere in $\Omega $, and $v\left(
0\right) =v_{0}$ almost everywhere in $\Gamma $. Note that in this setting, $%
v_{0}$ need not be the trace of $u_{0}$ at the boundary. Thus, in this
context equation (\ref{dyn}) is interpreted as an additional parabolic
equation, acting now on the boundary $\Gamma $.
\end{definition}

The following result is a direct consequence of results contained in \cite[%
Section 2]{GW}. The proof is based on the application of a Galerkin
approximation scheme which is not standard due to the nature of the boundary
conditions (see, also, \cite{CGGM}).

\begin{theorem}
\label{weaksol}Let the assumptions of (\ref{n1}), (\ref{n2}) be satisfied.
For any given initial data $U_{0}\in \mathbb{X}^{2},$ the problem (\ref{she}%
), (\ref{dyn}), (\ref{ini}) has a unique weak solution which depends
continuously on the initial data in a Lipschitz way. The following estimate
holds: 
\begin{align}
& \left\Vert U\left( t\right) \right\Vert _{\mathbb{X}^{2}}^{2}+\int%
\limits_{t}^{t+1}\left( \left\Vert U\left( s\right) \right\Vert _{\mathbb{V}%
_{1}}^{2}+\left\Vert u\left( s\right) \right\Vert _{L^{p}\left( \Omega
\right) }^{p}\right) ds  \label{diss} \\
& \leq c\left( \left\Vert U\left( 0\right) \right\Vert _{\mathbb{X}%
^{2}}^{2}\right) e^{-\rho t}+c\left( 1+\left\Vert g\right\Vert _{L^{2}\left(
\Omega \right) }^{2}\right) ,  \notag
\end{align}%
for all $t\geq 0$, where the positive constants $c$, $\rho $ are independent
of time and initial data.
\end{theorem}

As a consequence, problem (\ref{she}), (\ref{dyn}), (\ref{ini}) defines a
(nonlinear) continuous semigroup $\mathcal{S}_{t}$ acting on the phase-space 
$\mathbb{X}^{2}$,%
\begin{equation*}
\mathcal{S}_{t}:\mathbb{X}^{2}\rightarrow \mathbb{X}^{2},\text{ }t\geq 0,
\end{equation*}%
given by%
\begin{equation*}
\mathcal{S}_{t}U_{0}=U\left( t\right) .
\end{equation*}

\begin{theorem}
\label{attr}Let $f$ satisfy assumptions (\ref{n1}), (\ref{n2}), let $g\in
L^{\infty }\left( \Omega \right) $ and $\Gamma \in \mathcal{C}^{2}$. Then, $%
\left\{ \mathcal{S}_{t}\right\} $ possesses the connected global attractor $%
\mathcal{A}_{W},$ which is a bounded subset of $\mathbb{V}_{2}\cap \mathbb{X}%
^{\infty }$. As a consequence, the global attractor contains only strong
solutions.
\end{theorem}

\begin{proof}
The existence of an absorbing set in $\mathbb{V}_{1}\cap L^{p}\left( \Omega
\right) $ and, hence, the existence of the global attractor $\mathcal{A}%
_{W}\subset \mathbb{V}_{1}$ follows from \cite[Theorem 2.8 and Corollary 3.11%
]{GW}. We will now show that the attractor is bounded in $\mathbb{X}^{\infty
}$, and also in $\mathbb{V}_{2}$. All the calculations below are formal.
However, they can be rigorously justified by means of the approximation
procedures devised in \cite{GW} and \cite{GG0} (cf. \cite{CGGM} also). From
now on, $c$ will denote a positive constant that is independent of time and
initial data, which only depends on the other structural parameters of the
problem, that is, $\left\vert \Omega \right\vert ,$ $\left\vert \Gamma
\right\vert ,$ $\eta _{i},$ $\nu ,$ $b,$ $\left\Vert g\right\Vert _{\infty }$
and $n$. Such a constant may vary even from line to line.

\noindent \textbf{Step 1}. We will first establish the existence of a
bounded absorbing set in $\mathbb{X}^{\infty }$. First note that by (\ref%
{diss}), there is a constant $C_{0}>0,$ independent of time and initial
data, such that for any bounded subset $B$ of $\mathbb{X}^{2}$, $\exists $ $%
\tau =\tau \left( \left\Vert B\right\Vert _{\mathbb{X}^{2}}\right) >0$ with%
\begin{equation}
\sup_{t\geq \tau }\left\Vert U\left( t\right) \right\Vert _{\mathbb{X}%
^{2}}\leq C_{0}.  \label{diss2}
\end{equation}%
We shall now perform an Alikakos-Moser iteration argument. We multiply (\ref%
{she}) by $\left\vert u\right\vert ^{r_{k}-2}u,$ $r_{k}:=2^{k},$ $k\geq 1,$
and integrate over $\Omega $. We obtain%
\begin{eqnarray}
&&\frac{1}{r_{k}}\frac{d}{dt}\left\Vert u\right\Vert
_{r_{k}}^{r_{k}}+\left\langle f\left( u\right) ,\left\vert u\right\vert
^{r_{k}-2}u\right\rangle _{2}+\nu \int_{\Omega }\nabla u\cdot \nabla \left(
\left\vert u\right\vert ^{r_{k}-2}u\right) dx  \label{eqn2} \\
&=&\nu \int_{\Gamma }\partial _{\mathbf{n}}u\left\vert \psi \right\vert
^{r_{k}-2}\psi dS+\left\langle \lambda u+g,\left\vert u\right\vert
^{r_{k}-2}u\right\rangle _{2}.  \notag
\end{eqnarray}%
Similarly, we multiply (\ref{dyn}) by $\left\vert \psi \right\vert
^{r_{k}-2}\psi /b$ and integrate over $\Gamma $. We have%
\begin{equation}
\frac{1}{br_{k}}\frac{d}{dt}\left\Vert \psi \right\Vert _{r_{k},\Gamma
}^{r_{k}}+\nu \int_{\Gamma }\partial _{\mathbf{n}}u\left\vert \psi
\right\vert ^{r_{k}-2}\psi dS=0.  \label{eqn3}
\end{equation}%
Adding the equalities (\ref{eqn2}), (\ref{eqn3}), we deduce%
\begin{align}
& \frac{1}{r_{k}}\frac{d}{dt}\left( \left\Vert U\right\Vert _{\mathbb{X}%
^{r_{k}}}^{r_{k}}\right) +\left\langle f\left( u\right) ,\left\vert
u\right\vert ^{r_{k}-2}u\right\rangle _{2}+\nu \int_{\Omega }\nabla u\cdot
\nabla \left( \left\vert u\right\vert ^{r_{k}-2}u\right) dx  \label{eqn4} \\
& =\left\langle \lambda u+g,\left\vert u\right\vert ^{r_{k}-2}u\right\rangle
_{2}.  \notag
\end{align}%
A simple manipulation of the third integral in (\ref{eqn4}), and employing
assumption (\ref{n2}) on $f$, we readily get the estimate:%
\begin{align}
& \frac{d}{dt}\left( \left\Vert U\right\Vert _{\mathbb{X}^{r_{k}}}^{r_{k}}%
\right) +\eta _{1}r_{k}\left\Vert u\right\Vert _{r_{k}+p-2}^{r_{k}+p-2}+\nu
\left( 2^{k}-1\right) 2^{2-k}\left\Vert \nabla \left\vert u\right\vert
^{2^{k-1}}\right\Vert _{2}^{2}  \label{eqn5} \\
& \leq r_{k}\left\langle \lambda u+g+C_{f},\left\vert u\right\vert
^{r_{k}-2}u\right\rangle _{2}.  \notag
\end{align}%
Next, using the fact that $\left\vert y\right\vert ^{r_{k}-2}\leq \left\vert
y\right\vert ^{r_{k}}+1,$ for all $k\geq 1$ and $y\in \mathbb{R}$, we
estimate the last term on the right-hand side of (\ref{eqn5}),%
\begin{equation}
\left\langle \lambda u+g+C_{f},\left\vert u\right\vert
^{r_{k}-2}u\right\rangle _{2}\leq c\left( \left\Vert u\right\Vert
_{r_{k}}^{r_{k}}+1\right) ,  \label{eqn6}
\end{equation}%
for some positive constant $c$ that depends on $\lambda $ and the $L^{\infty
}$-norm of $g$, but is independent of $k$. On the other hand, it follows
from Gagliardo-Nirenberg inequality, and Young's inequality for $\varepsilon
\in \left( 0,1\right) ,$ that%
\begin{equation}
\left\Vert v\right\Vert _{2}\leq c\left\Vert v\right\Vert _{H^{1}\left(
\Omega \right) }^{n/\left( n+2\right) }\left\Vert v\right\Vert
_{1}^{1-n/\left( n+2\right) }\leq \varepsilon \left\Vert v\right\Vert
_{H^{1}\left( \Omega \right) }+c\varepsilon ^{-n/2}\left\Vert v\right\Vert
_{1},  \label{estt}
\end{equation}%
which implies%
\begin{equation*}
\left\Vert \nabla v\right\Vert _{2}^{2}\geq \frac{1-\varepsilon }{%
\varepsilon }\left\Vert v\right\Vert _{2}^{2}-c\varepsilon
^{-n/2-1}\left\Vert v\right\Vert _{1}^{2}.
\end{equation*}%
Note that the estimate (\ref{estt})\ remains valid if one replaces the $%
L^{2}\left( \Omega \right) $ and $L^{1}\left( \Omega \right) $-norms by the $%
L^{2}\left( \Gamma \right) $ and $L^{1}\left( \Gamma \right) $-norms,
respectively, and $n$ by $n-1$, respectively. Setting $v=\left\vert
u\right\vert ^{r_{k-1}}$ in the above inequality, noting that $\left(
2^{k}-1\right) 2^{2-k}\geq 2,$ for each $k,$ and the fact that $Tr_{D}$ maps 
$H^{1}\left( \Omega \right) $ boundedly into $L^{2}\left( \Gamma \right) $,
we can estimate the gradient term in (\ref{eqn5}) in terms of%
\begin{equation*}
c\frac{1-\varepsilon }{\varepsilon }\left( \left\Vert u\right\Vert
_{r_{k}}^{r_{k}}+\left\Vert \psi \right\Vert _{r_{k},\Gamma }^{r_{k}}\right)
-c\varepsilon ^{-n/2-1}\left( \left\Vert \left\vert u\right\vert
^{r_{k-1}}\right\Vert _{1}^{2}+\left\Vert \left\vert \psi \right\vert
^{r_{k-1}}\right\Vert _{1,\Gamma }^{2}\right) .
\end{equation*}%
(see, e.g., \cite[Chapter 5]{Ma}). This estimate together with (\ref{eqn5}),
(\ref{eqn6}) yield%
\begin{align}
& \frac{d}{dt}\left( \left\Vert U\right\Vert _{\mathbb{X}^{r_{k}}}^{r_{k}}%
\right) +c\left( \nu \frac{1-\varepsilon }{\varepsilon }-r_{k}\right) \left(
\left\Vert u\right\Vert _{r_{k}}^{r_{k}}+\left\Vert \psi \right\Vert
_{r_{k},\Gamma }^{r_{k}}\right)   \label{eqn7} \\
& \leq c\varepsilon ^{-n/2-1}\left( \left\Vert \left\vert u\right\vert
^{r_{k-1}}\right\Vert _{1}^{2}+\left\Vert \left\vert \psi \right\vert
^{r_{k-1}}\right\Vert _{1,\Gamma }^{2}\right) +cr_{k},  \notag
\end{align}%
for all $k\geq 1,$ where $c>0$ is independent of $k$.

We shall now make use of an iterative argument to deduce the existence of a
bounded absorbing set in $\mathbb{X}^{r_{k}},$ for all $k\geq 1$. Thus,
noting that $r_{k}\leq \left( r_{k}\right) ^{n/2+1}$, then choosing $%
\varepsilon =\delta /r_{k}$ with small $\delta =\delta \left( \nu \right) >0$
such that%
\begin{equation*}
\left( \nu \frac{1-\varepsilon }{\varepsilon }-r_{k}\right) \geq r_{k},
\end{equation*}%
and setting%
\begin{equation}
\mathcal{Y}_{k}\left( t\right) :=\int_{\Omega }\left\vert u\left( t,\cdot
\right) \right\vert ^{r_{k}}dx+\int_{\Gamma }\left\vert \psi \left( t,\cdot
\right) \right\vert ^{r_{k}}\frac{dS}{b}=\left\Vert U\right\Vert _{\mathbb{X}%
^{r_{k}}}^{r_{k}},  \label{def}
\end{equation}%
from (\ref{eqn7}) we derive the following estimate:%
\begin{equation}
\partial _{t}\mathcal{Y}_{k}\left( t\right) +cr_{k}\mathcal{Y}_{k}\left(
t\right) \leq c\left( r_{k}\right) ^{n/2+1}\left( \mathcal{Y}_{k-1}\left(
t\right) +1\right) ^{2}.  \label{e11}
\end{equation}%
Let us now take two positive constants $t,$ $\mu $ such that $t-\mu
/r_{k}>0, $ for all $k\geq 1$. Their precise values will be chosen later. We
claim that%
\begin{equation}
\mathcal{Y}_{k}\left( t\right) \leq M_{k}\left( t,\mu \right) :=c\left(
r_{k}\right) ^{n/2+1}(\sup_{s\geq t-\mu /r_{k}}\mathcal{Y}_{k-1}\left(
s\right) +1)^{2}  \label{claim}
\end{equation}%
holds for $\mathcal{Y}_{k},$ defined by (\ref{def}) and $k\geq 1$. To this
end, let $\zeta \left( s\right) $ be a positive function $\zeta :\mathbb{R}%
_{+}\rightarrow \left[ 0,1\right] $ such that $\zeta \left( s\right) =0$ for 
$s\in \left[ 0,t-\mu /r_{k}\right] ,$ $\zeta \left( s\right) =1$ if $s\in %
\left[ t,+\infty \right) $ and $\left\vert d\zeta /ds\right\vert \leq Cr_{k}$%
, if $s\in \left( t-\mu /r_{k},t\right) $. We define $Z_{k}\left( s\right)
=\zeta \left( s\right) \mathcal{Y}_{k}\left( s\right) $ and notice that%
\begin{equation*}
\frac{d}{ds}Z_{k}\left( s\right) \leq cr_{k}Z_{k}\left( s\right) +\zeta
\left( s\right) \frac{d}{ds}\mathcal{Y}_{k}\left( s\right) .
\end{equation*}%
Combining this estimate with (\ref{e11}) and noticing that $Z_{k}\leq 
\mathcal{Y}_{k}$, we deduce the following estimate for $Z_{k}$:%
\begin{equation}
\frac{d}{ds}Z_{k}\left( s\right) +cr_{k}Z_{k}\left( s\right) \leq
M_{k}\left( t,\mu \right) ,\text{ for all }s\in \left[ t-\mu /r_{k},+\infty
\right) .  \label{e12}
\end{equation}%
Integrating (\ref{e12}) with respect to $s$ from $t-\mu /r_{k}$ to $t$ and
taking into account the fact that $Z_{k}\left( t-\mu /r_{k}\right) =0,$ we
obtain that $\mathcal{Y}_{k}\left( t\right) =Z_{k}\left( t\right) \leq
M_{k}\left( t,\mu \right) \left( 1-e^{-C\mu }\right) $, which proves the
claim (\ref{claim}).

Let now $\tau ^{^{\prime }}>\tau >0$ be given with $\tau $ as in (\ref{diss2}%
), and define $\mu =2(\tau ^{^{\prime }}-\tau ),$ $t_{0}=\tau ^{^{\prime }}$
and $t_{k}=t_{k-1}-\mu /r_{k},$ $k\geq 1$. Using (\ref{claim}), we have%
\begin{equation}
\sup_{t\geq t_{k-1}}\mathcal{Y}_{k}\left( t\right) \leq c\left( r_{k}\right)
^{n/2+1}(\sup_{s\geq t_{k}}\mathcal{Y}_{k-1}\left( s\right) +1)^{2},\text{ }%
k\geq 1.  \label{e13}
\end{equation}%
Note that from (\ref{diss2}), we have $\left( \sup_{s\geq t_{1}=\tau }%
\mathcal{Y}_{1}\left( s\right) +1\right) \leq C_{0}+1=:\overline{C}$. Thus,
we can iterate in (\ref{e13}) with respect to $k\geq 1$ and obtain that%
\begin{align*}
\sup_{t\geq t_{k-1}}\mathcal{Y}_{k}\left( t\right) & \leq \left[ c\left(
r_{k}\right) ^{n/2+1}\right] \left[ c\left( r_{k-1}\right) ^{n/2+1}\right]
^{2}\cdot ...\cdot \left[ c\left( r_{1}\right) ^{n/2+1}\right] ^{2^{k}}(%
\overline{C})^{r_{k}} \\
& \leq c^{A_{k}}2^{B_{k}n/2+1}\left( \overline{C}\right) ^{r_{k}},
\end{align*}%
where%
\begin{equation}
A_{k}:=1+2+2^{2}+...+2^{k}\leq 2^{k}\sum_{i=1}^{\infty }\frac{1}{2^{i}}
\label{ak}
\end{equation}%
and%
\begin{equation}
B_{k}:=k+2\left( k-1\right) +2^{2}\left( k-2\right) +...+2^{k}\leq
2^{k}\sum_{i=1}^{\infty }\frac{i}{2^{i}}.  \label{bk}
\end{equation}%
Therefore,%
\begin{equation}
\sup_{t\geq t_{0}}\mathcal{Y}_{k}\left( t\right) \leq \sup_{t\geq t_{k-1}}%
\mathcal{Y}_{k}\left( t\right) \leq c^{A_{k}}2^{B_{k}\left( n/2+1\right) }%
\overline{C}^{r_{k}}.  \label{e14}
\end{equation}%
Since the series in (\ref{ak}) and\ (\ref{bk})\ are convergent, we can take
the $r_{k}$-root on both sides of (\ref{e14}) and let $k\rightarrow +\infty $%
. We deduce%
\begin{equation}
\sup_{t\geq t_{0}=\tau ^{^{\prime }}}\left\Vert U\left( t\right) \right\Vert
_{\mathbb{X}^{\infty }}\leq \lim_{k\rightarrow +\infty }\sup_{t\geq
t_{0}}\left( \mathcal{Y}_{k}\left( t\right) \right) ^{1/r_{k}}\leq C_{1},
\label{linf}
\end{equation}%
for some positive constant $C_{1}$ independent of $t,$ $k$, $U$ and initial
data.

\noindent \textbf{Step 2. }We claim that there is a positive constant $%
C_{2}, $ independent of time and initial data, and there exists $\tau
^{^{\prime \prime }}>0$ such that 
\begin{equation}
\left\Vert U\left( t\right) \right\Vert _{\mathbb{V}_{2}}\leq C_{2},\qquad 
\text{for all }\,t\geq \tau ^{^{\prime \prime }}.  \label{h2est}
\end{equation}%
Before we prove (\ref{h2est}), let us recall the following estimate (see 
\cite[Theorems 3.5, 3.10]{GW}):%
\begin{align}
& \sup_{t\geq \tau _{0}}\left( \left\Vert U\left( t\right) \right\Vert _{%
\mathbb{V}_{1}}^{2}+\left\Vert \partial _{t}u\left( t\right) \right\Vert
_{2}^{2}+\frac{1}{b}\left\Vert \partial _{t}\psi \left( t\right) \right\Vert
_{2,\Gamma }^{2}\right)  \label{rec} \\
& +\sup_{t\geq \tau _{0}}\int_{t}^{t+1}\left\Vert \partial _{t}u\left(
s\right) \right\Vert _{H^{1}\left( \Omega \right) }^{2}ds  \notag \\
& \leq C_{3},  \notag
\end{align}%
for some positive constant $C_{3}$ that is independent of time and the
initial data. In order to deduce (\ref{h2est}) from (\ref{rec}) and (\ref%
{linf}), we need to differentiate (\ref{she}) and (\ref{dyn}) with respect
to time. This yields%
\begin{equation}
\partial _{t}^{2}u=\nu \Delta \partial _{t}u-f^{^{\prime }}\left( u\right)
\partial _{t}u+\lambda \partial _{t}u,\text{ }\left( \partial _{t}^{2}\psi
+\nu b\partial _{{\mathbf{n}}}\left( \partial _{t}u\right) \right) _{\mid
\Gamma }=0.  \label{diff}
\end{equation}%
Then, we multiply the first equation of (\ref{diff}) by $\partial
_{t}^{2}u(t)$ and integrate over $\Omega ,$ using the boundary condition of (%
\ref{diff}). After standard transformations, we obtain 
\begin{align*}
& \frac{1}{2}\frac{d}{dt}\left( \left\Vert \nabla \partial _{t}u\left(
t\right) \right\Vert _{2}^{2}\right) +\left\Vert \partial _{t}^{2}u\left(
t\right) \right\Vert _{2}^{2}+\frac{1}{b}\left\Vert \partial _{t}^{2}\psi
\left( t\right) \right\Vert _{2,\Gamma }^{2} \\
& =-\left\langle \left( f^{^{\prime }}\left( u\left( t\right) \right)
-\lambda \right) \partial _{t}u\left( t\right) ,\partial _{t}^{2}u\left(
t\right) \right\rangle _{2}.
\end{align*}%
Using H\"{o}lder and Young inequalities, we have%
\begin{align*}
& \frac{d}{dt}\left( \left\Vert \nabla \partial _{t}u\left( t\right)
\right\Vert _{2}^{2}\right) +\left\Vert \partial _{t}^{2}u\left( t\right)
\right\Vert _{2}^{2}+\frac{2}{b}\left\Vert \partial _{t}^{2}\psi \left(
t\right) \right\Vert _{2,\Gamma }^{2} \\
& \leq c\left( \left\Vert f^{^{\prime }}\left( u\left( t\right) \right)
\partial _{t}u\left( t\right) \right\Vert _{2}^{2}+\left\Vert \partial
_{t}u\left( t\right) \right\Vert _{2}^{2}\right) \\
& \leq Q\left( \left\Vert u\left( t\right) \right\Vert _{\infty }\right)
\left\Vert \partial _{t}u\left( t\right) \right\Vert _{2}^{2},
\end{align*}%
for some positive nondecreasing function $Q$ that depends only on $f$ and $c$%
. This estimate yields, owing to (\ref{linf}), (\ref{rec}), 
\begin{equation*}
\frac{d}{dt}\left\Vert \nabla \partial _{t}u\left( t\right) \right\Vert
_{2}^{2}\leq c.
\end{equation*}%
Then, we can apply the so-called uniform Gronwall's lemma (see, e.g., \cite[%
Chapter III, Lemma 1.1]{T})\ to find a time $\tau _{1}\geq 1$, depending on $%
\tau _{0}$ and $\tau ,$ such that 
\begin{equation}
\left\Vert \nabla \partial _{t}u\left( t\right) \right\Vert _{2}^{2}\leq
c,\qquad \text{for all }\,t\geq \tau _{1}.  \label{4.29}
\end{equation}%
Therefore, (\ref{4.29}) and (\ref{rec}) allow us to deduce from (\ref{she})
and (\ref{dyn}), via standard elliptic regularity, the following estimate 
\begin{equation}
\left\Vert u\left( t\right) \right\Vert _{H^{2}\left( \Omega \right)
}^{2}\leq c,\qquad \forall \,t\geq \tau _{1}.  \label{4.30}
\end{equation}%
Summing up, we conclude by observing that (\ref{h2est}) follows from (\ref%
{4.30}) and the boundedness of the trace map $Tr_{D}:H^{2}\left( \Omega
\right) \rightarrow H^{3/2}\left( \Gamma \right) $. This completes the proof
of the theorem.
\end{proof}

\begin{remark}
The proof of Theorem \ref{attr} shows how to get an absorbing set in $%
\mathbb{V}_{2}$. Because of this, we can also prove the existence of a
global attractor for the dynamical system $\left( \left\{ \mathcal{S}%
_{t}\right\} _{t\geq 0},\mathbb{V}_{1}\right) .$
\end{remark}

\begin{theorem}
\label{reg2}If $\Omega $ is a bounded $\mathcal{C}^{\infty }$-domain, and $%
f,g$ are $\mathcal{C}^{\infty }$ functions, then the global attractor $%
\mathcal{A}_{W}$ is a bounded subset of $\mathbb{V}_{k},$ for every $k\geq 1$%
. In particular, if $U\in \mathcal{A}_{W}$ then $u\in \mathcal{C}^{\infty
}\left( \overline{\Omega }\right) .$
\end{theorem}

The proof of this result is standard and follows by successive time
differentiation of the equations in (\ref{diff}) and an induction argument.
We omit the details.

To prove the finite dimensionality of the global attractor $\mathcal{A}_{W}$%
, we can proceed in two different ways. One way is to establish the
existence of a more refined object called exponential attractor $\mathcal{E}%
_{W}$, whose existence proof is often based on proper forms of the so-called
squeezing/smoothing property for the differences of solutions. This can be
done by assuming smoother nonlinearities, i.e., $f\in C^{2}\left( \mathbb{R}%
\right) $ (see, e.g., \cite{GG0, GG1}). This has been carried out in \cite%
{GG0}, and references therein, for a system of reaction-diffusion equations
with dynamic boundary conditions of the form (\ref{dyn}), without relating
the attractor dimension to the physical parameters of the problem. However,
since we wish to find explicit estimates of fractal or/and Hausdorff
dimension of $\mathcal{A}_{W}$, we shall employ the classical machinery for
proving the finite dimensionality of the global attractor $\mathcal{A}_{W}.$
This is based on the so-called volume contraction arguments and requires the
associated solution semigroup $\mathcal{S}_{t}$ to be (uniformly quasi-)
differentiable with respect to the initial data, at least on the attractor
(see, e.g., \cite{BV}).

We give without proof the following result, which follows as a consequence
of the boundedness of $\mathcal{A}_{W}$ into $\mathbb{V}_{2}\cap \mathbb{X}%
^{\infty }.$

\begin{proposition}
\label{propd}Provided that $f\in C^{2}\left( \mathbb{R}\right) $ satisfies
the conditions (\ref{n1}) and (\ref{n2}), the flow $\mathcal{S}_{t}$
generated by the reaction-diffusion equation (\ref{she}) and dynamic
boundary condition (\ref{dyn}) is uniformly differentiable on $\mathcal{A}%
_{W},$ with differential%
\begin{equation}
\mathbf{L}\left( t,U\left( t\right) \right) :\Theta =\binom{\xi _{1}}{\xi
_{2}}\in \mathbb{X}^{2}\mapsto V=\binom{v}{\varphi }\in \mathbb{X}^{2},
\label{4.9}
\end{equation}%
where $V$ is the unique solution to%
\begin{align}
\partial _{t}v& =\nu \Delta v-f^{^{\prime }}\left( u\left( t\right) \right)
v+\lambda v,\text{ }\left( \partial _{t}\varphi +\nu b\partial _{\mathbf{n}%
}v\right) _{\mid \Gamma }=0,  \label{var} \\
V\left( 0\right) & =\Theta .  \notag
\end{align}%
Furthermore, $\mathbf{L}\left( t,U\left( t\right) \right) $ is compact for
all $t>0.$
\end{proposition}

The main result of this section is

\begin{theorem}
\label{main1}Let the assumptions of Proposition \ref{propd} be satisfied.
The fractal dimension of $\mathcal{A}_{W}$ admits the estimate%
\begin{equation}
\dim _{F}\mathcal{A}_{W}\leq c_{0}\left( 1+\frac{c_{f}+\lambda }{C_{W}\left(
\Omega ,\Gamma \right) \nu }\right) ^{n-1},\text{ for }n\geq 2  \label{updim}
\end{equation}%
and%
\begin{equation}
\dim _{F}\mathcal{A}_{W}\leq c_{0}\left( 1+\frac{c_{f}+\lambda }{C_{D}\left(
\Omega \right) \nu }\right) ^{1/2},\text{ for }n=1,  \label{updim2}
\end{equation}%
where $c_{0}$ depends on the shape of $\Omega $ only. The positive constants 
$C_{W},C_{D}$ depend only on $n,$ $\Omega ,$ $\Gamma $, $b$ and are given in
the Appendix.
\end{theorem}

\begin{proof}
In order to deduce (\ref{updim})-(\ref{updim2}), it is sufficient (see,
e.g., \cite[Chapter III, Definition 4.1]{CV})\ to estimate the $j$-trace of
the operator%
\begin{equation*}
\mathbf{L}\left( t,U\left( t\right) \right) =\left( 
\begin{array}{cc}
\nu \Delta -f^{^{\prime }}\left( u\left( t\right) \right) +\lambda I & 0 \\ 
-b\nu \partial _{\mathbf{n}} & 0%
\end{array}%
\right) .
\end{equation*}%
We have%
\begin{align*}
Trace\left( \mathbf{L}\left( t,U\left( t\right) \right) Q_{m}\right) &
=\sum_{j=1}^{m}\left\langle \mathbf{L}\left( t,U\left( t\right) \right)
\varphi _{j},\varphi _{j}\right\rangle _{\mathbb{X}^{2}} \\
& =\sum_{i=1}^{m}\left\langle \nu \Delta \varphi _{j},\varphi
_{j}\right\rangle _{2}-\sum_{i=1}^{m}\left\langle \nu \partial _{\mathbf{n}%
}\varphi _{j},\varphi _{j}\right\rangle _{2,\Gamma } \\
& -\sum_{i=1}^{m}\left\langle f^{^{\prime }}\left( u\left( t\right) \right)
\varphi _{j},\varphi _{j}\right\rangle _{2}+\sum_{i=1}^{m}\lambda
\left\langle \varphi _{j},\varphi _{j}\right\rangle _{2},
\end{align*}%
where the set of vector-valued functions $\varphi _{j}\in \mathbb{X}^{2}\cap 
\mathbb{V}_{1}$ is an orthonormal basis in $Q_{m}\mathbb{X}^{2}$. Since the
family $\varphi _{j}$ is orthonormal in $Q_{m}\mathbb{X}^{2},$ using
assumption (\ref{n1}) on $f$ (i.e., $f^{^{\prime }}\left( y\right) \geq
-c_{f},$ for all $y\in \mathbb{R}$), we find%
\begin{equation*}
Trace\left( \mathbf{L}\left( t,U\right) Q_{m}\right) \leq -\nu
\sum_{i=1}^{m}\left\Vert \nabla \varphi _{j}\right\Vert _{2}^{2}+\left(
c_{f}+\lambda \right) m.
\end{equation*}%
Let $n\geq 2$. From (\ref{LTineq}) (see Appendix, Proposition \ref{LT}), we
obtain%
\begin{align*}
Trace\left( \mathbf{L}\left( t,U\right) Q_{m}\right) & \leq -\nu
c_{1}C_{W}\left( \Omega ,\Gamma \right) m^{\frac{1}{n-1}+1}+\left( c_{1}\nu
C_{W}\left( \Omega ,\Gamma \right) +c_{f}+\lambda \right) m \\
& =:\rho \left( m\right) .
\end{align*}%
The function $\rho \left( y\right) $ is concave. The root of the equation $%
\rho \left( d\right) =0$ is%
\begin{equation*}
d^{\ast }=\left( 1+\frac{c_{f}+\lambda }{\nu c_{1}C_{W}\left( \Omega ,\Gamma
\right) }\right) ^{n-1}.
\end{equation*}%
Thus, we can apply \cite[Corollary 4.2 and Remark 4.1]{CV} to deduce that $%
\dim _{F}\mathcal{A}_{W}\leq d^{\ast },$ from which (\ref{updim}) follows.
The case $n=1$ is similar.
\end{proof}

\begin{remark}
Concerning the reaction-diffusion equation (\ref{she}), we can also handle
dynamic boundary conditions that involve surface diffusion: 
\begin{equation}
\partial _{t}u-\alpha \Delta _{\Gamma }u+b\nu \partial _{{\mathbf{n}}}\phi
=0,\text{ on }\Gamma ,  \label{3.32}
\end{equation}%
where $\alpha >0$ and $\Delta _{\Gamma }$ is the Laplace-Beltrami operator
on $\Gamma $. Our method of establishing upper bounds, comparable to the
bounds (\ref{updim})-(\ref{updim2}), for the dimension of the global
attractor can be also extended to this case as well. The details will appear
elsewhere.
\end{remark}

\section{Lower bounds on the dimension}

\label{lbb}

Lower bounds on the dimension of the global attractor are usually based on
the observation that the unstable manifold of any equilibrium of the system
is always contained in the global attractor (see, e.g., \cite{BV}). Thus, a
lower bound on the dimension of the attractor $\mathcal{A}_{W}$ can be found
by analyzing the dimension of an unstable manifold associated with a
constant equilibrium $Z$ for (\ref{she}), (\ref{dyn}). We begin by assuming
that $g$ is constant, for the sake of simplicity. Steady-state solutions of (%
\ref{she}), (\ref{dyn}) satisfy%
\begin{equation*}
L_{0}\left( u\right) :=\nu \Delta u-f\left( u\right) +\lambda u-g=0,\text{ }%
\left( \partial _{\mathbf{n}}u\right) _{\mid \Gamma }=0.
\end{equation*}%
We seek a solution of this system $U=\binom{u}{Tr_{D}\left( u\right) }\in 
\mathbb{X}^{2}$\ which coincides with a constant vector $Z=\mathbf{c}=\binom{%
c}{c},$ $c$ is a constant. Such a stationary solution satisfies the equation 
$\overline{L}_{0}\left( z\right) :=-f\left( z\right) +\lambda z-g=0.$ Since%
\begin{equation*}
f\left( y\right) y\geq \eta _{1}\left\vert y\right\vert ^{p}-C_{f},\text{
for }p>2,
\end{equation*}%
we have $\overline{L}_{0}\left( z\right) z\leq -\widetilde{\eta }%
_{1}\left\vert z\right\vert ^{p}+\widetilde{C}_{f},$ for some positive
constants $\widetilde{\eta }_{1},\widetilde{C}_{f}$. Thus, $\overline{L}%
_{0}\left( z\right) z<0$ on the interval $I_{R}=\left( -R,R\right) ,$ if $R$
is large enough. It follows that $\overline{L}_{0}\left( z\right) =0$ has at
least one solution $Z=Z\left( \lambda \right) $ (see, e.g., \cite[Chapter III%
]{CV}). By the implicit function theorem, this solution is of order
\thinspace $1/\lambda $ for sufficiently large $\lambda .$

Now fix this solution. In order to find a lower bound on the dimension of
the global attractor $\mathcal{A}_{W},$ it suffices to establish a lower
bound for $\dim E_{+}\left( Z\right) ,$ where $E_{+}\left( Z\right) $ is an
invariant subspace of $\mathbf{L}\left( Z\right) ,$ which corresponds to%
\begin{equation*}
\mathbf{L}\left( Z\right) W=\binom{\nu \Delta w-f^{^{\prime }}\left(
z\right) w+\lambda w}{-b\nu \partial _{\mathbf{n}}w}
\end{equation*}%
with $\sigma \left( \mathbf{L}\left( Z\right) \right) \subset \left\{ \zeta
:\zeta >0\right\} $. We note that $\left( \mathbf{L}\left( Z\right) ,D\left( 
\mathbf{L}\left( Z\right) \right) \right) $\ is self-adjoint on $X^{2}$\
with spectrum contained in $\left( -\infty ,c_{f}+\lambda \right] .$

The main result of this section is the following.

\begin{theorem}
\label{lbbb}Let $f\in C^{2}\left( \mathbb{R}\right) $ satisfy assumptions (%
\ref{n1})-(\ref{n2}). There exist a positive constant $c_{0}$, depending on $%
f,$ $g$ and the shape of $\Omega ,$ independent of $\lambda ,$ $\nu ,$ $b$, $%
\left\vert \Omega \right\vert ,$ $\left\vert \Gamma \right\vert ,$ such that%
\begin{equation*}
\dim _{F}\mathcal{A}_{W}\geq \dim _{H}\mathcal{A}_{W}\geq \dim E_{+}\left(
Z\right) \geq c_{0}\left( \frac{\lambda }{C_{W}\left( \Omega ,\Gamma \right)
\nu }\right) ^{n-1},
\end{equation*}%
for $n\geq 2$. In one space dimension, the same estimate is valid with $%
C_{W} $ replaced by $C_{D}$ and $n-1,$ replaced by $1/2$, respectively.
\end{theorem}

\begin{proof}
For a fixed constant solution $Z=\mathbf{c}$ of $\overline{L}_{0}\left(
z\right) =0$ and sufficiently large $\lambda \geq 1,$ we have $\chi \left(
\lambda \right) :=-f^{^{\prime }}\left( z\right) +\lambda >0$.

Let $\left\{ \varphi _{j}\left( x\right) \right\} _{ji\in \mathbb{N}_{0}}$
be an orthonormal basis in $\mathbb{X}^{2}$ consisting of eigenfunctions of
the Wentzell Laplacian $\Delta _{W}$ (see Appendix, Theorem \ref{sanalysis2}%
),%
\begin{equation}
\Delta _{W}\varphi _{j}=\Lambda _{j}\varphi _{j},\text{ }j\in \mathbb{N}_{0},%
\text{ }\varphi _{j}\in D\left( \Delta _{W}\right) \cap C\left( \overline{%
\Omega }\right)  \label{evseq}
\end{equation}%
such that%
\begin{equation*}
0=\Lambda _{0}<\Lambda _{1}\leq \Lambda _{2}\leq ...\leq \Lambda _{,j}\leq
\Lambda _{j+1}\leq ....
\end{equation*}%
We shall seek for eigenvectors $W_{j}=\binom{w_{j}}{Tr_{D}\left(
w_{j}\right) }\in \mathbb{X}^{2}$, of the form $w_{j}\left( x\right)
=\varphi _{j}\left( x\right) p_{j},$ $p_{j}\in \mathbb{R}$, satisfying
equation%
\begin{equation}
\mathbf{L}\left( Z\right) W_{j}=\zeta _{j}W_{j},\text{ }W_{j}\in D\left( 
\mathbf{L}\left( Z\right) \right) :=D\left( \Delta _{W}\right) .  \label{eee}
\end{equation}%
Note that for $W_{j}\in D\left( \mathbf{L}\left( Z\right) \right) \subset 
\mathbb{V}_{1},$ the trace of $w_{j}$ makes sense as an element of $%
H^{1/2}\left( \Gamma \right) $. Substituting such $w_{j}$ into (\ref{eee}),
taking into account (\ref{evseq}) and the fact that%
\begin{equation*}
\mathbf{L}\left( Z\right) W_{j}=-\nu \Delta _{W}W_{j}+\Pi _{\lambda }W_{j},%
\text{ }\Pi _{\lambda }W_{j}:=\binom{\chi \left( \lambda \right) w_{j}}{0},
\end{equation*}%
we obtain the equation%
\begin{equation}
\left( -\nu \Lambda _{j}I+\Pi _{\lambda }\right) p_{j}=\zeta _{j}p_{j},\text{
}\Pi _{\lambda }=\left( 
\begin{array}{cc}
\chi \left( \lambda \right) & 0 \\ 
0 & 0%
\end{array}%
\right) .  \label{pj}
\end{equation}%
A nonzero $p_{j}$ exists if $\zeta =\zeta _{j}$ is a root of the equation%
\begin{equation}
\det \left( -\nu \Lambda _{j}I+\Pi _{\lambda }-\zeta I\right) =0,\text{ }%
\zeta >0.  \label{dett}
\end{equation}%
When $\nu =0,$ this equation has at least one root $\zeta >0$ since $\chi
=\chi \left( \lambda \right) >0$ (in fact, $\zeta =\chi \left( \lambda
\right) $). Therefore, there exists $\delta >0$ such that when $\nu \Lambda
_{j}<\delta ,$ the equation (\ref{dett}) has a root $\zeta _{j}=\zeta
_{j}\left( \nu \right) $ with $\zeta _{j}>0$. Therefore, to any such root $%
\zeta _{j}$, we can assign a nontrivial $p_{j},$ which is a solution of (\ref%
{pj}), and thus an eigenvector $W_{j}=\binom{w_{j}}{Tr_{D}w_{j}},$ $%
w_{j}=\varphi _{j}p_{j}$. Let us now compute how many $j$'s satisfy the
inequality $\nu \Lambda _{j}<\delta $. The asymptotic behavior of $\Lambda
_{j}$ is $\Lambda _{j}\sim C_{W}\left( \Omega ,\Gamma \right) j^{1/\left(
n-1\right) }$ as $j\rightarrow \infty $ (see, Appendix, Theorem \ref%
{asymptotic2}). The last inequality certainly holds when%
\begin{equation*}
1\leq j\leq c_{1}\delta ^{n-1}\left( C_{W}\nu \right) ^{1-n}=c_{2}\left( 
\frac{1}{C_{W}\nu }\right) ^{n-1},\text{ for }n\geq 2
\end{equation*}%
and%
\begin{equation*}
1\leq j\leq c_{1}\delta ^{1/2}\left( C_{D}\nu \right) ^{-1/2}=c_{2}\left( 
\frac{1}{C_{D}\nu }\right) ^{1/2},\text{ for }n=1.
\end{equation*}%
The positive constants $c_{1},$ $c_{2}$ depend on $\lambda .$ In order to
get more explicit estimates for $c_{1},$ $c_{2}$, it is left to remark that
equation (\ref{dett}) may be rewritten in the form%
\begin{equation*}
\det \left( -\nu \Lambda _{j}\lambda ^{-1}I+\lambda ^{-1}\Pi _{\lambda
}-\zeta _{1}I\right) =0
\end{equation*}%
with $\zeta _{1}=\lambda ^{-1}\zeta ,$ and to observe that a solution of
this equation clearly exists if $\nu \Lambda _{j}\lambda ^{-1}\leq \delta ,$
for sufficiently large $\lambda $ and small $\delta $. Employing the
asymptotic formula for $\Lambda _{j}$ once again, we find 
\begin{equation*}
1\leq j\leq c_{1}^{^{\prime }}\delta ^{n-1}\lambda ^{n-1}\left( C_{W}\nu
\right) ^{1-n}=c_{2}^{^{\prime }}\left( \frac{\lambda }{C_{W}\nu }\right)
^{n-1},\text{ for }n\geq 2
\end{equation*}%
and%
\begin{equation*}
1\leq j\leq c_{1}^{^{\prime }}\delta ^{1/2}\lambda ^{1/2}\left( C_{D}\nu
\right) ^{-1/2}=c_{2}^{^{\prime }}\left( \frac{\lambda }{C_{D}\nu }\right)
^{1/2},\text{ for }n=1.
\end{equation*}%
It follows that 
\begin{equation*}
\dim E_{+}\left( Z\left( \lambda \right) \right) \geq c_{2}^{^{\prime
}}\left( \frac{\lambda }{C_{W}\nu }\right) ^{n-1},\text{ for }n\geq 2
\end{equation*}%
and 
\begin{equation*}
\dim E_{+}\left( Z\left( \lambda \right) \right) \geq c_{2}^{^{\prime
}}\lambda ^{1/2}\left( C_{D}\nu \right) ^{-1/2},
\end{equation*}%
in one space dimension. The proof is complete.
\end{proof}

\section{Concluding remarks}

\label{cr}

In the textbook literature on theoretical geophysics, it was traditional to
use a Robin boundary condition with a nonlinear heat equation to describe
temperature variations at the upper surface of the ocean \cite{LMT, LMT2}.
But it was recognized that this was not always the physically correct
boundary condition \cite{MW2}. Among its applicability to a wide range of
phenomena, including phase-transitions in fluids, and so on \cite{GG0, QWS},
the reaction-diffusion system (\ref{she})-(\ref{dyn}) has important
applications in climatology and is essentially used to determine large and
rapid temperature changes in the ocean's surface as a response to changes
into deep water formations \cite{MW2}. In this paper, we provide explicit
bounds for the dimension of the attractor for this system and study the
effect of the dynamic term $b^{-1}\partial _{t}u$, representing change in
thermal energy\ in an infinitesimal layer near the surface. Unlike the
previous results, the dimension of the attractor is proportional to the
surface area $\left\vert \Gamma \right\vert ,$ for large domains $\Omega $
and fixed parameters $\nu ,$ $\lambda $ and $b$. Moreover, all the constants
involved in our estimates are given in an explicit form. We also observe
that in the case without $b^{-1}\partial _{t}u$ in (\ref{dyn}), \thinspace
i.e., $b=+\infty ,$ the dimension of the attractor is much larger (and
proportional to the volume $\left\vert \Omega \right\vert $ of $\Omega $)
than the dimension of the global attractor for the same system when $0<b\neq
+\infty $. Thus, we observe that the addition of the dynamic term $%
b^{-1}\partial _{t}u$, $b>0$ drastically changes the situation. This is a
remarkable fact that can have a profound effect onto the long-term dynamics
of other systems that are subject to dynamic boundary conditions of this
form. We will investigate these effects for other systems, such as the B\'{e}%
nard problem for nonlinear heat conduction, in forthcoming papers. Finally,
we note that it is also possible to extend the results of this paper to the
case when the boundary $\Gamma $\ consists of two disjoint open subsets $%
\Gamma _{1}$\ and $\Gamma _{2}$, each $\overline{\Gamma }_{i}\backprime
\Gamma _{i}$\ is a $S$-null subset of $\Gamma $\ and $\Gamma =\overline{%
\Gamma }_{1}\cup \overline{\Gamma }_{2}$\ with $\Gamma _{1}\subsetneqq
\Gamma $, such that $u$\ satisfies a Dirichlet boundary condition on $\Gamma
_{1}$\ and a dynamic boundary condition on $\Gamma _{2}$. We will come back
to this issue in a forthcoming article.

\section{Appendix}

\label{ap}

In this section, we shall recall several important results concerning a
certain realization of $L=\nu \Delta $ with the Wentzell boundary condition (%
\ref{wbc}). We have the following.

\begin{theorem}
\label{Wentzell2}Let $\Omega $ be a bounded open set of $\mathbb{R}^{n}$
with Lipschitz boundary $\Gamma $. Assume that $b>0$ and $0\leq q\in
L^{\infty }\left( \Omega \right) $. Define the operator $\Delta _{W}$ on $%
\mathbb{X}^{2},$ by%
\begin{equation}
\Delta _{W}\binom{u_{1}}{u_{2}}:=\binom{-\Delta u_{1}+q\left( x\right) u_{1}%
}{b\partial _{\mathbf{n}}u_{1}},  \label{A_Wentzell1}
\end{equation}%
with%
\begin{equation}
D\left( \Delta _{W}\right) :=\left\{ U=\binom{u_{1}}{u_{2}}\in \mathbb{V}%
_{1}:-\Delta u_{1}\in L^{2}\left( \Omega \right) ,\text{ }\partial _{\mathbf{%
n}}u_{1}\in L^{2}\left( \Gamma ,\frac{dS}{b}\right) \right\} .
\label{A_Wentzell2}
\end{equation}%
Then, $\left( \Delta _{W},D\left( \Delta _{W}\right) \right) $ is
self-adjoint on $\mathbb{X}^{2}.$ Moreover, the resolvent operator $\left(
I+\Delta _{W}\right) ^{-1}\in \mathcal{L}\left( \mathbb{X}^{2}\right) $ is
compact.
\end{theorem}

We refer the reader to \cite{CFGGOR, GG_b, GGGRW} for an extensive survey of
recent results concerning the "Wentzell" Laplacian $\Delta _{W}$.

The eigenvalue problem associated with the operator $\Delta _{W}$ is given
by $\Delta _{W}\varphi =\Lambda \varphi ;$ this leads to the following
spectral problem for the perturbed Laplacian%
\begin{equation}
-\Delta \varphi +q\left( x\right) \varphi =\Lambda \varphi \text{ in }\Omega
,  \label{sp1}
\end{equation}%
with a boundary condition that depends on the eigenvalue $\Lambda $
explicitly:%
\begin{equation}
b\partial _{\mathbf{n}}\varphi =\Lambda \varphi \text{ on }\Gamma .
\label{sp2}
\end{equation}%
Such a function $\varphi $ will be called an eigenfunction associated with $%
\Lambda $ and the set of all eigenvalues $\Lambda $ of (\ref{sp1})-(\ref{sp2}%
) will be denoted by $\Lambda _{W}.$ Let $\varphi _{j}$ and $\Lambda _{W,j}$%
, $j\in J$, denote all the eigenfunctions and eigenvalues of (\ref{sp1})-(%
\ref{sp2}). We have the following (see, e.g., \cite{BBR, VVi}).

\begin{theorem}
\label{sanalysis1}Let $q\geq 0$ with $\int\limits_{\Omega }q\left( x\right)
dx>0$. Then, there exists a sequence of numbers%
\begin{equation}
0<\Lambda _{W,1}\leq \Lambda _{W,2}\leq ...\leq \Lambda _{W,j}\leq \Lambda
_{W,j+1}\leq ...,  \label{seq}
\end{equation}%
converging to $+\infty $, with the following properties:

(a) The spectrum of $\Delta _{W}$ is given by%
\begin{equation*}
\sigma \left( \Delta _{W}\right) =\left\{ \Lambda _{W,j}\right\} _{j\in 
\mathbb{N}},
\end{equation*}%
and each number $\Lambda _{W,j},$ $j\in \mathbb{N},$ is an eigenvalue for $%
\Delta _{W}$ of finite multiplicity.

(b) There exists a countable family of orthonormal eigenfunctions for $%
\Delta _{W}$ which spans $\mathbb{X}^{2}$. More precisely, there exists a
collection of functions $\left\{ \varphi _{j}\right\} _{j\in \mathbb{N}}$
with the following properties:%
\begin{align}
\varphi _{j}& \in D\left( \Delta _{W}\right) \text{ and }\Delta _{W}\varphi
_{j}=\Lambda _{W,j}\varphi _{j},\text{ }j\in \mathbb{N},  \label{prop} \\
\left\langle \varphi _{j},\varphi _{k}\right\rangle _{\mathbb{X}^{2}}&
=\delta _{jk}\text{, }j,k\in \mathbb{N}\text{,}  \notag \\
\mathbb{X}^{2}& =\oplus \overline{lin.span\left\{ \varphi _{j}\right\}
_{j\in \mathbb{N}}}\text{ (orthogonal direct sum).}  \notag
\end{align}

(c) If $\Gamma $ is Lipschitz, then every eigenfunction $\varphi _{j}\in 
\mathbb{V}_{1}$, and in fact $\varphi _{j}\in C(\overline{\Omega })\cap
C^{\infty }(\Omega )$, for every $j$. If $\Gamma $ is of class $C^{2}$, then
every eigenfunction $\varphi _{j}\in \mathbb{V}_{1}\cap C^{2}\left( 
\overline{\Omega }\right) ,$ for every $j.$

(d) The following min-max principle holds:%
\begin{equation}
\Lambda _{W,j}=\min_{\substack{ Y_{j}\subset \mathbb{V}_{1},  \\ \dim
Y_{j}=j }}\max_{0\neq \varphi \in Y_{j}}R_{W}\left( \varphi ,\varphi \right)
,\text{ }j\in \mathbb{N}\text{,}  \label{minmax}
\end{equation}%
where the Rayleigh quotient $R_{W}$, for the perturbed Wentzell operators,
is given by%
\begin{equation}
R_{W}\left( \varphi ,\varphi \right) :=\frac{\left\Vert \nabla \varphi
\right\Vert _{2}^{2}+\left\langle q\left( x\right) \varphi ,\varphi
\right\rangle _{2}}{\left\Vert \varphi \right\Vert _{\mathbb{X}^{2}}^{2}},%
\text{ }0\neq \varphi \in \mathbb{V}_{1}.  \label{rqw}
\end{equation}
\end{theorem}

Concerning the case $q\equiv 0$, we have the following.

\begin{theorem}
\label{sanalysis2}Let $q\equiv 0.$ Then, there exists a sequence of numbers%
\begin{equation*}
0=\Lambda _{W,0}<\Lambda _{W,1}\leq \Lambda _{W,2}\leq ...\leq \Lambda
_{W,j}\leq \Lambda _{W,j+1}\leq ...,
\end{equation*}%
converging to $+\infty $, with the following properties:

(a) The spectrum of $\Delta _{W}$ is given by%
\begin{equation*}
\sigma \left( \Delta _{W}\right) =\left\{ \Lambda _{W,j}\right\} _{j\in 
\mathbb{N\cup }\left\{ 0\right\} },
\end{equation*}%
and each number $\Lambda _{W,j},$ $j\in \mathbb{N}_{0}=\mathbb{N\cup }%
\left\{ 0\right\} ,$ is an eigenvalue for $\Delta _{W}$ of finite
multiplicity. The eigenvalue $\Lambda _{W,0}$ is simple and its associated
eigenfunction is of constant sign.

(b) There exists a countable family of orthonormal eigenfunctions for $%
\Delta _{W}$ which spans $\mathbb{X}^{2}$. More precisely, the same
conclusion (b) of Theorem \ref{sanalysis1} holds in this case as well.
Finally, both conclusions (c) and (d) in Theorem \ref{sanalysis1} hold in
the case $q\equiv 0$ as well.
\end{theorem}

The asymptotic behavior of the eigenvalues $\Lambda _{W,j},$ as $%
j\rightarrow \infty ,$ was established in \cite{Fran, Fran2}. We refer the
reader to \cite{GG_b} for more details about the Wentzell Laplacian and
other generalizations. Let $J=\mathbb{N}_{0}$ or $\mathbb{N}$, according to
whether $q=0$ or $q>0$ respectively. Set%
\begin{equation*}
C_{D}\left( \Omega \right) :=\frac{\left( 2\pi \right) ^{2}}{\left(
v_{n}\left\vert \Omega \right\vert \right) ^{2/n}}\text{ and }C_{S}\left(
\Gamma \right) =\frac{2\pi }{\left( v_{n-1}\left\vert \Gamma \right\vert
\right) ^{1/\left( n-1\right) }}.
\end{equation*}%
Here $v_{n}$ denotes the volume of the unit ball in $\mathbb{R}^{n}$, and we
recall that $\left\vert \Omega \right\vert $ stands for the $n$-dimensional
Euclidean volume of $\Omega $, while $\left\vert \Gamma \right\vert $ stands
for the usual $\left( n-1\right) $-dimensional Lebesgue surface measure on $%
\Gamma $.

We summarize these results in the following.

\begin{theorem}
\label{asymptotic2}The eigenvalue sequence $\left\{ \Lambda _{W,j}\right\}
_{j\in J}$ of the (un)perturbed Wentzell Laplacian $\Delta _{W}$ satisfies:

(i) For $n\geq 2$, we have%
\begin{equation}
\Lambda _{W,j}=C_{W}\left( \Omega ,\Gamma \right) j^{1/\left( n-1\right)
}+o\left( j^{1/\left( n-1\right) }\right) ,\text{ as }j\rightarrow +\infty ,
\label{ev}
\end{equation}%
for some%
\begin{equation}
C_{W}\left( \Omega ,\Gamma \right) \in \left\{ 
\begin{array}{cc}
bC_{S}\left( \Gamma \right) \left[ 2^{-1/\left( n-1\right) },1\right] , & 
\text{for }n\geq 3, \\ 
\left[ \frac{C_{D}\left( \Omega \right) C_{S}\left( \Gamma \right) }{2\left(
b^{-1}C_{D}\left( \Omega \right) +C_{S}\left( \Gamma \right) \right) },\min
\left\{ C_{D}\left( \Omega \right) ,bC_{S}\left( \Gamma \right) \right\} %
\right] , & \text{for }n=2.%
\end{array}%
\right.  \label{constant}
\end{equation}

(ii) For $n=1,$ we have%
\begin{equation}
\Lambda _{W,j}=C_{D}\left( \Omega \right) j^{2}+o\left( j^{2}\right) ,\text{
as }j\rightarrow +\infty .  \label{ev2}
\end{equation}
\end{theorem}

The following version of the Lieb--Thirring inequality is essential.

\begin{proposition}
\label{LT}Let $\omega _{j},$ $1\leq j\leq m,$ be a finite family of $\mathbb{%
V}_{1},$ which is orthonormal in $\mathbb{X}^{2}$. We have%
\begin{equation}
\sum_{i=1}^{m}\left\Vert \nabla \omega _{j}\right\Vert _{2}^{2}\geq
c_{1}C_{W}\left( \Omega ,\Gamma \right) \left( m^{\frac{1}{n-1}+1}-m\right) .
\label{LTineq}
\end{equation}%
The constant $c_{1}>0$ depends only on $n$ and the shape of $\Omega ,$ and
is independent of the size of $\Omega ,$ $\Gamma ,$ of $m,$ and of the $%
\omega _{j}$'s.
\end{proposition}

\begin{proof}
Let $B_{W}:=\Delta _{W}+C_{W}\left( \Omega ,\Gamma \right) I$ and let $%
D\left( B_{W}\right) =D\left( \Delta _{W}\right) $. By Theorems \ref%
{Wentzell2}, \ref{sanalysis1}, $B_{W}$ is a linear positive unbounded
self-adjoint operator on $\mathbb{X}^{2},$ such that $B_{W}^{-1}$ is
compact. Thus, we can apply the abstract result of \cite[Chapter VI, Lemma
2.1]{T} to deduce that%
\begin{align}
\sum_{i=1}^{m}\left( \left\Vert \nabla \omega _{j}\right\Vert
_{2}^{2}+C_{W}\left\Vert \omega _{j}\right\Vert _{\mathbb{X}^{2}}^{2}\right)
& =\sum_{i=1}^{m}\left\langle B_{W}\omega _{j},\omega _{j}\right\rangle _{%
\mathbb{X}^{2}}  \label{LTpp} \\
& \geq \Lambda _{W,1}\left( B_{W}\right) +\Lambda _{W,2}\left( B_{W}\right)
+...+\Lambda _{W,m}\left( B_{W}\right)  \notag \\
& \geq C_{W}\left( 1^{1/\left( n-1\right) }+2^{1/\left( n-1\right)
}+...+m^{1/\left( n-1\right) }\right)  \notag \\
& \geq c_{0}C_{W}m^{\frac{1}{n-1}+1},  \notag
\end{align}%
since, by (\ref{ev})-(\ref{ev2}), $\Lambda _{W,j}\left( B_{W}\right) \geq
C_{W}\left( \Omega ,\Gamma \right) j^{1/\left( n-1\right) },$ for all $j,$
and some positive constant $c_{0}$ (indeed, we have $\Lambda _{W,j}\left(
B_{W}\right) =\Lambda _{W,j}\left( \Delta _{W}\right) +C_{W}$)$.$ Thus, the
proof of (\ref{LTineq}) follows immediately from (\ref{LTpp}).
\end{proof}

\end{document}